# D'Alembert's solution of fractional wave equations using complex fractional transformation

**Uttam Ghosh[1a], Md Ramjan Ali[1b], Santanu Raut[2], Susmita Sarkar[1c] and Shantanu Das[3]**

[1]Department of Applied Mathematics, University of Calcutta, Kolkata, India
[1a]email: uttam_math@yahoo.co.in
[1b]email: ramjan.azad@gmail.com
[1c]email: susmita62@yahoo.co.in
[2]Mathabhanga College, Cooch Behar, West Bengal, India
email: raut_santanu@yahoo.com
[3]Reactor Control Systems Design Section E & I Group BARC Mumbai India
email: shantanu@barc.gov.in

## Abstract

Fractional wave equation arises in different type of physical problems such as the vibrating strings, propagation of electro-magnetic waves, and for many other systems. The exact analytical solution of the fractional differential equation is difficult to find. Usually Laplace-Fourier transformation method, along with methods where solutions are represented in series form is used to find the solution of the fractional wave equation. In this paper we describe the D'Alembert's solution of the fractional wave equation with the help of complex fractional transform method. We demonstrate that using this fractional complex transformation method, we obtain the solutions easily as compared to fractional method of characteristics; and get the solution in analytical form. We show that the solution to the fractional wave equation manifests as travelling waves with scaled coordinates, depending on the considered fractional order value.

## 1.0 Introduction

Fractional calculus is one of the oldest inventions of modern mathematics. It originated from the letter of Leibnitz to L'Hospital in 1695 [1-2, 30]. Fractional differential models are used to describe the physical phenomena those have memory or hereditary effect [2] such as visco-elasticity, anomalous diffusion, fluid mechanics, biology, acoustics, di-electric relaxations, dynamics of super-capacitors , control theory etc [2,4].

The most commonly used definition of fractional derivative is the Riemann-Liouville (R-L) fractional derivative [2, 5-10]. The differential equation formed with R-L type fractional derivative has initial/ boundary conditions defined via fractional order derivative [2]. In real life problem it is difficult to formulate the initial/ boundary conditions with fractional derivative. Therefore, the researchers are using the Caputo [2, 11] or Jumarie [12] type fractional derivative to formulate the problems. To solve the fractional differential equation usually methods used are the Variation Iteration Method [13-14], exp-method [15], fractional sub-equation method [16-18], solution using Mittag-Leffler functions as eigen-function [2, 19, and 30] for fractional differential operator, Complex fractional transform method [21-25]. Among all these methods the complex fractional transformation method is one of the simplest methods to solve the Fractional Partial Differential Equation (FPDE). This transform transforms the fractional partial differential equation to a classical partial differential equation and thus the solution procedure become simple.

The D'Alembert's solution of Cauchy Problem of the wave equation plays an important role in different type of vibration and wave propagation problems. Usually Laplace-Fourier transform method gives the exact solution of the wave equation, in closed form by use of Mittag-Leffler and M-Wright functions, where obtaining inverse Laplace-Fourier transforms are difficult using contour integration or Berberan-Santo's method [30]. The other techniques to get the solution of fractional differential equation are by the methods such as Variation Iteration Method [13-14], Adomian Decomposition Method [2] etc; where the solutions obtained are in form of series solution; (which is analytical approximate of closed form solution). In this paper we find the D'Alembert's solution of fractional wave equation using the complex transformation method as exact analytic solution in closed form.

One of the motivations for using fractional calculus in physical systems is that the space and time variables which we often deal exhibit coarse-grained phenomena. This means infinitesimal quantities cannot be arbitrarily taken to zero-rather they are non-zero with a minimum spread. This type of non-zero spread arises in the microscopic to mesoscopic levels of system dynamics. This means that if we denote $x$ as the point in space and $t$ as the point in time, then limit of the differentials $dx$ (and $dt$) cannot be taken zero. To take this concept of coarse graining into account, use the infinitesimal quantities as $(\Delta x)^\alpha$ (and $(\Delta t)^\alpha$) with $0 < \alpha < 1$; called as 'fractional differentials'. For arbitrarily small $\Delta x$ and $\Delta t$ (tending towards zero), these 'fractional' differentials are greater than $\Delta x$ (and $\Delta t$) i.e. $(\Delta x)^\alpha > \Delta x$ and $(\Delta t)^\alpha > \Delta t$. This way of defining the fractional differentials helps us to use fractional derivatives in the study of dynamic systems.

The coordinate $x$ in $x$-space, originates from the differential $dx$, as its integration i.e. $\int_0^x d\xi = x$. Now with a differential $(dx)^\alpha$, with $0 < \alpha < 1$, we have $(dx)^\alpha > dx$, where $\int_0^x (d\xi)^\alpha \sim x^\alpha$ [30, 31]. That is the space is transformed to a fractal space where the coordinate $x$ is transformed to $x^\alpha$. This is coarse graining phenomena in particular scale of observation. Here we come across the fractal space time, where the normal classical differentials $dx$ and $dt$, cannot be taken arbitrarily to zero. Thus in these cases the concept of classical differentiability is lost [30]. The fractional order $\alpha$ is related to roughness character of the space-time i.e. the fractal dimension [30, 31]. We will show that solution of fractional wave equations have component of travelling wave in the scaled coordinates $\sim x^\alpha$ and $t^\alpha$.

Organization of the paper is as follows: in section 2.0 we describe some review of the fractional derivative and complex fractional transformation. The section 3.0 describes the method to find solution of $\alpha$ − th and $2\alpha$ − th order wave equation using complex fractional transformation. In section 4.0 the existence and uniqueness of the solution is described. Finally conclusion is given in section 5.0; followed by references.

## 2.0 Review of fractional calculus and Complex transform method

In this section we describe some basic definitions of the fractional derivatives and the theory of complex fractional transformation.

a) **Fractional derivative**

The most commonly used fractional derivative is the Riemann–Liouville (R-L) [2, 4, 30], that is defined in integral form as

$$^{RL}_{t_0}D^\alpha_t f(t) = \frac{1}{\Gamma(n-\alpha)} \frac{d^n}{dt^n} \int_{t_0}^{t} (t-\xi)^{n-\alpha-1} f(\xi) d\xi \tag{2.1}$$

$$t > t_0, \quad n-1 < \alpha \leq n$$

In terms of this definition fractional derivative of the constant $K$ is $^{RL}_{t_0}D^\alpha_t K = \frac{K}{\Gamma(1-\alpha)}(t-t_0)^{-\alpha} \neq 0$ [2, 4, and 30]. In classical calculus the fractional derivative of the constant equal to zero. To overcome this difference M. Caputo [11] proposed the fractional order derivative. Let $\alpha$ be a positive real number and $n$ be a positive integer satisfying $n-1 < \alpha \leq n$. Let $f^{(n)}(t)$ exists (that is $n$−th order ordinary derivative exists). Then $\alpha$-th order fractional derivative in Caputo sense is defined by

$$^{C}_{t_0}D^\alpha_t f(t) = ^{RL}_{t_0}I^{n-\alpha}_t \left[ f^{(n)}(t) \right] = \frac{1}{\Gamma(n-\alpha)} \int_{t_0}^{t} (t-\xi)^{n-\alpha-1} f^{(n)}(\xi) d\xi, \tag{2.2}$$

$$t > t_0, \quad n-1 < \alpha \leq n$$

where $^{RL}_{t_0}I^\alpha_t$ is the R-L integral operator[2, 4, 30], defined as follows

$$^{RL}_{t_0}I^\alpha_t \left[ f(t) \right] = \frac{1}{\Gamma(\alpha)} \int_{t_0}^{t} (t-\xi)^{\alpha-1} f(\xi) d\xi, \quad \alpha > 0$$

In terms of the Caputo derivative fractional derivative of constant ($K$) is zero i.e. $^{C}_{a}D^\alpha_t K = 0$. But this Caputo definition is applicable only for differentiable functions [2, 4, 30]. Jumarie [12, 30] proposed the modified R-L derivative of continuous (but not necessarily differentiable) function $f(x)$ in the range $0 \leq x \leq a$ in the following form, (with assumption that $f(0)$ is finite)

$$_0^J D_x^\alpha [f(x)] = f^{(\alpha)}(x) = \begin{cases} \dfrac{1}{\Gamma(-\alpha)} \displaystyle\int_0^x (x-\xi)^{-\alpha-1} f(\xi) d\xi, & \alpha < 0 \\ \dfrac{1}{\Gamma(1-\alpha)} \dfrac{d}{dx} \displaystyle\int_0^x (x-\xi)^{-\alpha} \big(f(\xi) - f(0)\big) d\xi, & 0 < \alpha < 1 \\ \big(f^{(\alpha-n)}(x)\big)^{(n)}, & n \leq \alpha < n+1, \quad n \geq 1 \end{cases} \quad (2.3)$$

On the other hand the general form of fractional differential equation of two independent variables $x$, $t$ and one dependent variable $u(x,t)$ is defined as function $\Phi$ as follows

$$\Phi\big(u, D_x^\alpha u, D_t^\alpha u, \ldots, D_x^{2\alpha} u, D_t^{2\alpha} u\big) = 0, \quad 0 < \alpha < 1 \quad (2.4)$$

where $D_x^\alpha u = \frac{\partial^\alpha u(x,t)}{\partial x^\alpha}$ denotes the Jumarie fractional partial derivative[22] of the form,

$$\frac{\partial^\alpha u(x,t)}{\partial x^\alpha} = \frac{1}{\Gamma(1-\alpha)} \frac{d}{dx} \int_0^x (x-\xi)^{-\alpha} \big(u(\xi,t) - u(0,t)\big) d\xi, \quad 0 < \alpha \leq 1 \quad (2.5)$$

The function $u(x,t)$ in (2.5) is continuous (but not necessarily differentiable) function. The Jumarie fractional derivative of a constant is zero which overcomes the discomfort of RL derivative.

### b) The fractional Complex Transformation

Let the differential equation be of the form, with $D^\alpha$ as Jumarie fractional derivative operator

$$D^\alpha u = f(x), \qquad D^\alpha = \frac{d^\alpha}{dx^\alpha} \quad (2.6)$$

The complex fractional transformation for this differential equation (2.6) is defined as, $X = \frac{(px)^\alpha}{\Gamma(1+\alpha)}$ where $p$ is constant [22]. From Fractional Taylor series (of Jumarrie type [30]) we have a conversion formula that is $\alpha! dx \simeq x^{(\alpha)}(t)(dt)^\alpha$. From this have following conversion formula

$$\alpha! dx \simeq \left(\frac{d^\alpha x}{(dt)^\alpha}\right)(dt)^\alpha = d^\alpha x; \qquad \alpha! = \Gamma(\alpha+1)$$

$$d^\alpha x = \big(\Gamma(\alpha+1)\big) dx$$

The relation $d^\alpha x \simeq \big(\Gamma(\alpha+1)\big) dx$, as we get from above, is the conversion formula for fractional differential $d^\alpha f$ to $df$. Using this Jumarie conversion formula $d^\alpha u = \alpha! du = \Gamma(1+\alpha) du$ and then doing obvious manipulations as depicted in following steps, and thereafter using the formula $\frac{d^\alpha}{dx^\alpha} ax^\alpha = a(\Gamma(1+\alpha))$ (comes from Euler formula [2, 30] i.e. $\frac{d^\alpha}{dx^\alpha} x^\beta = \frac{\Gamma(\beta+1)}{\Gamma(\beta+1-\alpha)} x^{\beta-\alpha}$) we get

$$\frac{d^\alpha u}{dx^\alpha} = \frac{\Gamma(1+\alpha)du}{dx^\alpha}$$

$$= \frac{du}{dX} \frac{\Gamma(1+\alpha)dX}{dx^\alpha}$$

$$= \left(\frac{du}{dX}\right)\left(\frac{d^\alpha X}{dx^\alpha}\right)$$

$$= \left(\frac{du}{dX}\right)\left(\frac{d^\alpha \left(\frac{px^\alpha}{\Gamma(1+\alpha)}\right)}{dx^\alpha}\right)$$

$$= p^\alpha \frac{du}{dX}$$

With this above derivation the equation (2.6) i.e. $D_x^\alpha u = f(x)$ reduces to the integer order ordinary differential equation, that is

$$p^\alpha \frac{du}{dX} = F(X), \quad X = \frac{(px)^\alpha}{\Gamma(1+\alpha)} .$$

For the fractional partial differential of following form

$$D_t^\alpha u(x,t) + D_x^\alpha u(x,t) = f(x,t) \qquad (2.7)$$

where $0 < \alpha \leq 1$ and $D_t^\alpha u(x,t) = \frac{\partial^\alpha u}{\partial t^\alpha}, D_x^\alpha u(x,t) = \frac{\partial^\alpha u}{\partial x^\alpha}$; by using the complex fractional transformations i.e. $X = \frac{(px)^\alpha}{\Gamma(1+\alpha)}, T = \frac{(qt)^\alpha}{\Gamma(1+\alpha)}$ we get the following transformed expressions for fractional partial derivative operators

$$\frac{\partial^\alpha u}{\partial x^\alpha} = \frac{\Gamma(1+\alpha)\partial u}{\partial x^\alpha} = \frac{\partial u}{\partial X} \frac{\Gamma(1+\alpha)\partial X}{\partial x^\alpha} = p^\alpha \frac{\partial u}{\partial X}; \qquad \frac{\partial^\alpha u}{\partial t^\alpha} = q^\alpha \frac{\partial u}{\partial T}$$

Using the above transformation equation (2.7) i.e. $D_t^\alpha u(x,t) + D_x^\alpha u(x,t) = f(x,t)$; reduces to the integer order (classical) partial differential equation that is following

$$q^\alpha \frac{\partial u}{\partial T} + p^\alpha \frac{\partial u}{\partial X} = F(X,T)$$

## 3.0 Application of Complex fractional transformation

In the next two sub-sections we shall describe the solution of $\alpha-$th and $2\alpha-$th order wave equation using the complex fractional transformation.

### 3.1 Solution of $\alpha-$th order wave equation

Now consider the $\alpha-$th order fractional differential equation $D_t^\alpha u(x,t) + c^\alpha D_x^\alpha u(x,t) = 0$ that is of the following form

$$\frac{\partial^{\alpha} u(x,t)}{\partial t^{\alpha}} + c^{\alpha} \frac{\partial^{\alpha} u(x,t)}{\partial x^{\alpha}} = 0, \qquad 0 < \alpha \leq 1 \tag{3.1}$$

The initial condition of (3.1) is specified as; $u(x,0) = f\left(\frac{x^{\alpha}}{\Gamma(1+\alpha)}\right)$.

First we will use fractional method of characteristic [28, 29]. Using this we have for FPDE (3.1) the characteristic that are following

$$\left. \begin{aligned} \frac{du}{ds} &= 0 \\ \frac{(dt)^{\alpha}}{(\Gamma(1+\alpha))ds} &= 1 \\ \frac{(dx)^{\alpha}}{(\Gamma(1+\alpha))ds} &= c^{\alpha} \end{aligned} \right\} \tag{3.2}$$

Now we use the formula of integration w.r.t. $(dx)^{\alpha}$ [28, 29, and 30] as indicated below for integration of (3.2)

$$\int_0^x f(\xi)(d\xi)^{\alpha} \stackrel{\text{def}}{=} \alpha \int_0^x (x-\xi)^{\alpha-1} f(\xi) d\xi, \qquad 0 < \alpha < 1$$
$$= \frac{\alpha(\Gamma(\alpha))}{\Gamma(\alpha)} \int_0^x (x-\xi)^{\alpha-1} f(\xi) d\xi = \Gamma(\alpha+1)\left( {}_0 I_x^{\alpha}[f(x)] \right)$$

Where ${}_0 I_x^{\alpha}[f(x)] = \frac{1}{\Gamma(\alpha)} \int_0^x (x-\xi)^{\alpha-1} f(\xi) d\xi$ is Riemann-Liouville fractional integration of order $0 < \alpha < 1$. On making $f(x) = 1$ we get $\int_0^x (d\xi)^{\alpha} = x^{\alpha}$, $0 < \alpha \leq 1$.

Solving the first equation we get $u(x,t) = c_1$ (the reason that Jumarie Fractional Derivative of constant is zero) and from the last two equations we get [28, 29]

$$(dx)^{\alpha} = c^{\alpha} (dt)^{\alpha} \tag{3.3}$$

Integrating the expression (3.3) with the formula of integration w.r.t. $(dx)^{\alpha}$ [28, 29, 30] we get

$$\frac{x^{\alpha}}{\Gamma(1+\alpha)} - c^{\alpha} \frac{t^{\alpha}}{\Gamma(1+\alpha)} = c_2 \tag{3.4}$$

Hence the solution of the equation (3.1) can be written as

$$u(x,t) = \varphi\left( \frac{1}{\Gamma(1+\alpha)} x^{\alpha} - \frac{1}{\Gamma(1+\alpha)} c^{\alpha} t^{\alpha} \right) \tag{3.5}$$

where $\varphi(x)$ is an arbitrary function of single variable [28, 29]. Using the initial condition we get $\varphi(x) = f(x)$. Hence the solution of (3.1) is given by,

$$u(x,t) = f\left(\tfrac{1}{\Gamma(1+\alpha)} x^\alpha - \tfrac{1}{\Gamma(1+\alpha)} c^\alpha t^\alpha\right) \tag{3.6}$$

Now use the fractional complex transformation, with $\frac{p^\alpha x^\alpha}{\Gamma(1+\alpha)} = X$ and $\frac{q^\alpha t^\alpha}{\Gamma(1+\alpha)} = T$ developed by He [22-23]; and write the following expressions for fractional differential operators (that get changed to classical differential operator), as follows

$$\frac{\partial^\alpha u}{\partial t^\alpha} = q^\alpha \frac{\partial u}{\partial T}, \qquad \frac{\partial^\alpha u}{\partial x^\alpha} = p^\alpha \frac{\partial u}{\partial X} \tag{3.7}$$

With above change by use of fractional complex transformation, we obtain following equivalence in operator form

$$\frac{\partial^\alpha}{\partial t^\alpha} + c^\alpha \frac{\partial^\alpha}{\partial x^\alpha} \equiv q^\alpha \frac{\partial}{\partial T} + c^\alpha p^\alpha \frac{\partial}{\partial X}$$

Using the transformation (3.7) equation (3.1) reduces to the first order partial differential equation in the following form

$$\begin{aligned}
&\frac{\partial^\alpha u(x,t)}{\partial t^\alpha} + c^\alpha \frac{\partial^\alpha u(x,t)}{\partial x^\alpha} = 0, \qquad X = \frac{p^\alpha x^\alpha}{\Gamma(1+\alpha)}, \qquad T = \frac{q^\alpha t^\alpha}{\Gamma(1+\alpha)} \\
&q^\alpha \frac{\partial u(X,T)}{\partial T} + c^\alpha p^\alpha \frac{\partial u(X,T)}{\partial X} = 0
\end{aligned} \tag{3.8}$$

With choice of $p = q = 1$ (as these constants are scale factors) we get following

$$X = \frac{x^\alpha}{\Gamma(1+\alpha)} \qquad T = \frac{t^\alpha}{\Gamma(1+\alpha)} \qquad \frac{\partial u(X,T)}{\partial T} + c^\alpha \frac{\partial u(X,T)}{\partial X} = 0$$

with the initial condition $u(X,0) = f(X)$. The characteristics of (3.8) with $p = q = 1$ are following that is given by,

$$\frac{du}{ds} = 0; \qquad \frac{dX}{ds} = c^\alpha; \qquad \frac{dT}{ds} = 1$$

Solution of the first equation is $u(X,T) = c_1$ and from the last two equations after integration we get $X - c^\alpha T = c_2$, where $c_1$ and $c_2$ are constants. Hence solution of the partial differential equation (3.8) is $u(X,T) = \varphi(X - c^\alpha T)$. Using the initial conditions we get, $\varphi(X) = f(X)$. Hence the solution of (3.8) is given by

$$u(X,T) = f\left(X - c^\alpha T\right) \tag{3.9}$$

Recalling the original variable in (3.9) i.e. $\frac{x^\alpha}{\Gamma(1+\alpha)} = X$ and $\frac{t^\alpha}{\Gamma(1+\alpha)} = T$ we get solution of the $\alpha-$th order partial differential equation (3.1) is

$$u(X,T) = u(x,t) = f\left(\tfrac{1}{\Gamma(1+\alpha)}x^\alpha - \tfrac{1}{\Gamma(1+\alpha)}c^\alpha t^\alpha\right) \tag{3.10}$$

The solution (3.6) and solution (3.10) are the same expressions therefore the solution obtained by both the methods is same. We observe the solution obtained by using complex transform method is easy comparative to the earlier method i.e. fractional method of characteristics. The original initial function $f(x)$, travels with constant velocity in this case is $c^\alpha$, as it travels in scaled coordinates $\frac{1}{\Gamma(1+\alpha)}x^\alpha$ and $\frac{1}{\Gamma(1+\alpha)}t^\alpha$. These are fractional order travelling waves [31].

### 3.2 D'Alembert's Solution of $2\alpha-$th order fractional Wave Equation
Here we consider the one dimensional Fractional wave equation in the form,

$$\frac{\partial^{2\alpha} u(x,t)}{\partial t^{2\alpha}} - c^{2\alpha}\frac{\partial^{2\alpha} u(x,t)}{\partial x^{2\alpha}} = 0, \qquad 0 < \alpha \leq 1 \tag{3.11}$$

with the initial conditions expressed as

$$u(x,0) = f\left(\tfrac{x^\alpha}{\Gamma(1+\alpha)}\right); \quad \left.\frac{\partial^\alpha}{\partial t^\alpha}u(x,t)\right|_{t=0} = g\left(\tfrac{x^\alpha}{\Gamma(1+\alpha)}\right) \tag{3.11a}$$

where $c$ in (3.11) is a positive constant. This $c$ is physically interpreted as the wave speed. For example, if displacement of infinite string is represented by $u(x,t)$ then $c$ becomes the quantity $\sigma/\rho$ where $\sigma$ is the tension of the string and $\rho$ is the density. Now, the equation (3.11) can be written (in factored form) as following

$$\left(\frac{\partial^\alpha}{\partial t^\alpha} - c^\alpha \frac{\partial^\alpha}{\partial x^\alpha}\right)\left(\frac{\partial^\alpha}{\partial t^\alpha} + c^\alpha \frac{\partial^\alpha}{\partial x^\alpha}\right)u(x,t) = 0 \tag{3.12}$$

Let us consider the fractional complex transformation developed by He [22-23],

$$\frac{(px)^\alpha}{\Gamma(1+\alpha)} = X, \qquad \frac{(qt)^\alpha}{\Gamma(1+\alpha)} = T \tag{3.13}$$

Using this transformation we get,

$$\frac{\partial^\alpha u}{\partial t^\alpha} = q^\alpha \frac{\partial u}{\partial T}, \qquad \frac{\partial^\alpha u}{\partial x^\alpha} = p^\alpha \frac{\partial u}{\partial X}.$$

Implying the following operator equivalence

$$\frac{\partial^\alpha}{\partial t^\alpha} - c^\alpha \frac{\partial^\alpha}{\partial x^\alpha} \equiv q^\alpha \frac{\partial}{\partial T} - c^\alpha p^\alpha \frac{\partial}{\partial X}; \qquad \frac{\partial^\alpha}{\partial t^\alpha} + c^\alpha \frac{\partial^\alpha}{\partial x^\alpha} \equiv q^\alpha \frac{\partial}{\partial T} + c^\alpha p^\alpha \frac{\partial}{\partial X}$$

Using the transformation (3.13) the equation (3.12) reduces to the following form,

$$\left( q^\alpha \frac{\partial}{\partial T} - c^\alpha p^\alpha \frac{\partial}{\partial X} \right)\left( q^\alpha \frac{\partial}{\partial T} + c^\alpha p^\alpha \frac{\partial}{\partial X} \right) u(X,T) = 0 \qquad (3.14)$$

Under this transformation the boundary conditions reduces to following sets

$$u(X,0) = f(X); \qquad \frac{\partial}{\partial T} u(X,0) = \frac{g(X)}{q^\alpha} \qquad (3.15)$$

Let us consider the one factored term of (3.14) that is the following

$$\left( q^\alpha \frac{\partial}{\partial T} + c^\alpha p^\alpha \frac{\partial}{\partial X} \right) u(X,T) = v(X,T) \qquad (3.16)$$

then equation (3.14) reduces to

$$\left( q^\alpha \frac{\partial}{\partial T} - c^\alpha p^\alpha \frac{\partial}{\partial X} \right) v(X,T) = 0 \qquad (3.17)$$

Using the method of characteristic [26] the solution of the equation (3.17) can be written in the following form,

$$v(X,T) = f_1\left( q^\alpha X + c^\alpha p^\alpha T \right) \qquad (3.18)$$

Therefore from (3.16) $u(X,T)$ satisfies the following equation,

$$\left( q^\alpha \frac{\partial}{\partial T} + c^\alpha p^\alpha \frac{\partial}{\partial X} \right) u(X,T) = f_1\left( q^\alpha X + c^\alpha p^\alpha T \right) \qquad (3.19)$$

On the other hand the equation (3.14) can be rewritten as,

$$\left[ \left( q^\alpha \frac{\partial}{\partial T} + c^\alpha p^\alpha \frac{\partial}{\partial X} \right)\left( q^\alpha \frac{\partial}{\partial T} - c^\alpha p^\alpha \frac{\partial}{\partial X} \right) \right] u(X,T) = 0 \qquad (3.20)$$

Again consider $\left( q^\alpha \frac{\partial}{\partial T} - c^\alpha p^\alpha \frac{\partial}{\partial X} \right) u(X,T) = w(X,T)$ then equation (3.20) reduces to following

$$\left(q^\alpha \frac{\partial}{\partial T} + c^\alpha p^\alpha \frac{\partial}{\partial X}\right) w(X,T) = 0 \tag{3.21}$$

Thus (on similar lines as described above) solution of the equation (3.21) can be written in the following form

$$w(X,T) = g_1\left(q^\alpha X - c^\alpha p^\alpha T\right) \tag{3.22}$$

Therefore $u(X,T)$ satisfies the equation,

$$\left(q^\alpha \frac{\partial}{\partial T} - c^\alpha p^\alpha \frac{\partial}{\partial X}\right) u(X,T) = g_1\left(q^\alpha X - c^\alpha p^\alpha T\right) \tag{3.23}$$

where $f_1$ and $g_1$ are functions of single variable. Now adding the equations (3.19) and (3.23) we get following expression

$$2q^\alpha \frac{\partial}{\partial t} u(X,T) = f_1\left(q^\alpha X + c^\alpha p^\alpha T\right) + g_1\left(q^\alpha X - c^\alpha p^\alpha T\right) \tag{3.24}$$

Integrating (3.24) we can write the solution in the following form

$$u(X,T) = F\left(q^\alpha X + c^\alpha p^\alpha T\right) + G\left(q^\alpha X - c^\alpha p^\alpha T\right) \tag{3.25}$$

where $F$ and $G$ arbitrary functions of the single variable. Differentiating (3.25) with respect to $T$ we get

$$u_T(X,T) = c^\alpha p^\alpha F'\left(q^\alpha X + c^\alpha p^\alpha T\right) - c^\alpha p^\alpha G'\left(q^\alpha X - c^\alpha p^\alpha T\right) \tag{3.26}$$

Then from (3.25) and (3.26) and using initial conditions we get the following set of expressions

$$\left.\begin{array}{l} \dfrac{\partial}{\partial T} u(X,0) = c^\alpha p^\alpha \left(F'(X) - G'(X)\right) = \dfrac{g(X)}{q^\alpha} \\ u(X,0) = F(X) + G(X) = f(X) \\ F'(X) + G'(X) = f'(X) \end{array}\right\} \tag{3.27}$$

Solving the equations in (3.27) we get,

$$F'(X) = \tfrac{1}{2}\left(f'(X) - \tfrac{1}{c^\alpha p^\alpha q^\alpha} g(X)\right), \qquad G'(X) = \tfrac{1}{2}\left(f'(X) + \tfrac{1}{c^\alpha p^\alpha q^\alpha} g(X)\right) \tag{3.28}$$

Integrating (3.28) on both side we get following expressions

$$F(X) = \tfrac{1}{2}\left(f(X) - \tfrac{1}{c^\alpha p^\alpha q^\alpha}\int_{X_0}^{X} g(\xi)d\xi\right)$$
$$G(X) = \tfrac{1}{2}\left(f(X) + \tfrac{1}{c^\alpha p^\alpha q^\alpha}\int_{X_0}^{X} g(\xi)d\xi\right) \tag{3.29}$$

Thus from (3.25) we get the general solution of (3.1) which is

$$u(X,T) = \frac{1}{2}\left(f(q^\alpha X + c^\alpha p^\alpha T) + f(q^\alpha X - c^\alpha p^\alpha T)\right) + \frac{1}{2c^\alpha p^\alpha q^\alpha}\int_{q^\alpha X - c^\alpha p^\alpha T}^{q^\alpha X + c^\alpha p^\alpha T} g(\xi)d\xi. \tag{3.30}$$

Now transforming to the original variables i.e. $\frac{p^\alpha x^\alpha}{\Gamma(1+\alpha)}(\equiv X)$ and $\frac{q^\alpha t^\alpha}{\Gamma(1+\alpha)}(\equiv T)$ the solution (3.30) reduces to the following form,

$$u(x,t) = \frac{1}{2}\left[f\left(p^\alpha q^\alpha\left(\tfrac{1}{\Gamma(1+\alpha)}x^\alpha + \tfrac{1}{\Gamma(1+\alpha)}c^\alpha t^\alpha\right)\right) + f\left(p^\alpha q^\alpha\left(\tfrac{1}{\Gamma(1+\alpha)}x^\alpha - \tfrac{1}{\Gamma(1+\alpha)}c^\alpha t^\alpha\right)\right)\right]$$
$$+ \frac{1}{2c^\alpha p^\alpha q^\alpha}\int_{p^\alpha q^\alpha\left(\tfrac{1}{\Gamma(1+\alpha)}x^\alpha - \tfrac{1}{\Gamma(1+\alpha)}c^\alpha t^\alpha\right)}^{p^\alpha q^\alpha\left(\tfrac{1}{\Gamma(1+\alpha)}x^\alpha + \tfrac{1}{\Gamma(1+\alpha)}c^\alpha t^\alpha\right)} g(\xi)d\xi \tag{3.31}$$

With $p = q = 1$ (as these constants are scale factors) and $\alpha = 1$ the initial condition of (3.11) is $u(x,t)\big|_{t=0} = f(x)$. The equation (3.11) is thus a classical wave equation i.e. $\frac{\partial^2}{\partial t^2}u(x,t) - c^2\frac{\partial^2}{\partial x^2}u(x,t) = 0$. Assume the second initial condition as $\frac{\partial}{\partial t}u(x,t)\big|_{t=0} = 0 = g(x)$. Placing $p = q = 1$ and $g(x) = 0$ in (3.31), we get the following

$$u(x,t) = \tfrac{1}{2}f(x+ct) + \tfrac{1}{2}f(x-ct)$$

The above expression represents the initial condition (function) getting split into two with amplitude half each and travelling in $-x$ and $+x$ directions respectively with constant velocity $c$ in opposite directions. This is standard classical Cauchy problem of a classical wave equation. The solution (3.31) is D'Alembert solution of the Cauchy problem for one dimensional fractional wave equation given by (3.31). Now with $g(x) = 0$ in (3.31), and with $p = q = 1$ for solution to fractional wave equation (3.31), we get the two waves (defined by initial condition) travelling in opposite directions with velocity $c^\alpha$ in the scaled coordinates $\frac{1}{\Gamma(1+\alpha)}x^\alpha$ and $\frac{1}{\Gamma(1+\alpha)}t^\alpha$.

### 4.0 Existence and uniqueness

**Theorem:** Show that the problem (3.11) and (3.11a) in the domain $0 < x < \infty$  $0 \leq t \leq \mathbf{T}$ for a fixed $\mathbf{T} > 0$ well-posed for $f \in C^{2\alpha}(\mathbb{R})$ and $g \in C^\alpha(\mathbb{R})$.

**Proof:** The existence and uniqueness of the solutions of (3.11) with (3.11a) is directed from D'Alembert's solutions. Now from smoothness conditions $f \in C^{2\alpha}(\mathbb{R})$ and $g \in C^\alpha(\mathbb{R})$ implies that $u(x,t) \in C^{2\alpha}(\mathbb{R} \times (0,\infty)) \cap C^{\alpha}(\mathbb{R} \times (0,\infty))$. Thus like the classical integer order solution the solution (3.31) is generalized solution.

Now we shall prove the stability of the solution i.e. for a small change of the initial conditions give rise to small change of the solution. Let $u_i(x,t)$ be two solution of the problem (3.11) with two initial conditions $f_i, g_i$ for $i = 1,2$.

Now since $f_i \in C^{2\alpha}(\mathbb{R})$ and $g_i \in C^{\alpha}(\mathbb{R})$ for $i = 1,2$ with $\varepsilon > 0, \delta > 0$ and for all $x \in \mathbb{R}^+$
$|f_1(x) - f_2(x)| < \delta$; $|g_1(x) - g_2(x)| < \delta$. Then for all $x \in \mathbb{R}^+$ and $0 \leq t \leq T$ we have the following from (3.31)

$$|u_1(x,t) - u_2(x,t)| \leq \frac{\left|f_1\left(p^{\alpha}q^{\alpha}\left(\frac{1}{\Gamma(1+\alpha)}x^{\alpha} + \frac{1}{\Gamma(1+\alpha)}c^{\alpha}t^{\alpha}\right)\right) - f_2\left(p^{\alpha}q^{\alpha}\left(\frac{1}{\Gamma(1+\alpha)}x^{\alpha} + \frac{1}{\Gamma(1+\alpha)}c^{\alpha}t^{\alpha}\right)\right)\right|}{2}$$

$$+ \frac{\left|f_1\left(p^{\alpha}q^{\alpha}\left(\frac{1}{\Gamma(1+\alpha)}x^{\alpha} - \frac{1}{\Gamma(1+\alpha)}c^{\alpha}t^{\alpha}\right)\right) - f_2\left(p^{\alpha}q^{\alpha}\left(\frac{1}{\Gamma(1+\alpha)}x^{\alpha} - \frac{1}{\Gamma(1+\alpha)}c^{\alpha}t^{\alpha}\right)\right)\right|}{2}$$

$$+ \frac{1}{2c^{\alpha}p^{\alpha}q^{\alpha}} \int_{p^{\alpha}q^{\alpha}\left(\frac{1}{\Gamma(1+\alpha)}x^{\alpha} - \frac{1}{\Gamma(1+\alpha)}c^{\alpha}t^{\alpha}\right)}^{p^{\alpha}q^{\alpha}\left(\frac{1}{\Gamma(1+\alpha)}x^{\alpha} + \frac{1}{\Gamma(1+\alpha)}c^{\alpha}t^{\alpha}\right)} |g_1(\xi) - g_2(\xi)| d\xi$$

$$< \delta(1 + T^{\alpha}) < \varepsilon$$

if $\delta < \frac{1}{(1+T^{\alpha})}\varepsilon$. Therefore for all $x \in \mathbb{R}^+$ and $0 \leq t \leq T$ we have $|u_1(x,t) - u_2(x,t)| < \varepsilon$. Hence the theorem is proved.

**Example**: consider the fractional wave equation in the form

$$\frac{\partial^{2\alpha}u(x,t)}{\partial t^{2\alpha}} - c^{2\alpha}\frac{\partial^{2\alpha}u(x,t)}{\partial x^{2\alpha}} = 0, \qquad 0 < \alpha \leq 1 \qquad (3.32)$$

with the initial conditions

(1) $u(x,0) = f(x) = x^2$ and $D_t^{\alpha}u(x,t)\big|_{t=0} = g(x) = \sin(x)$

for first case represented; the solution is represented by figure-1; and with initial conditions

(2) $u(x,0) = f(x) = 0$ and $D_t^{\alpha}u(x,t)\big|_{t=0} = g(x) = \sin x$

the solution represented by figure-2

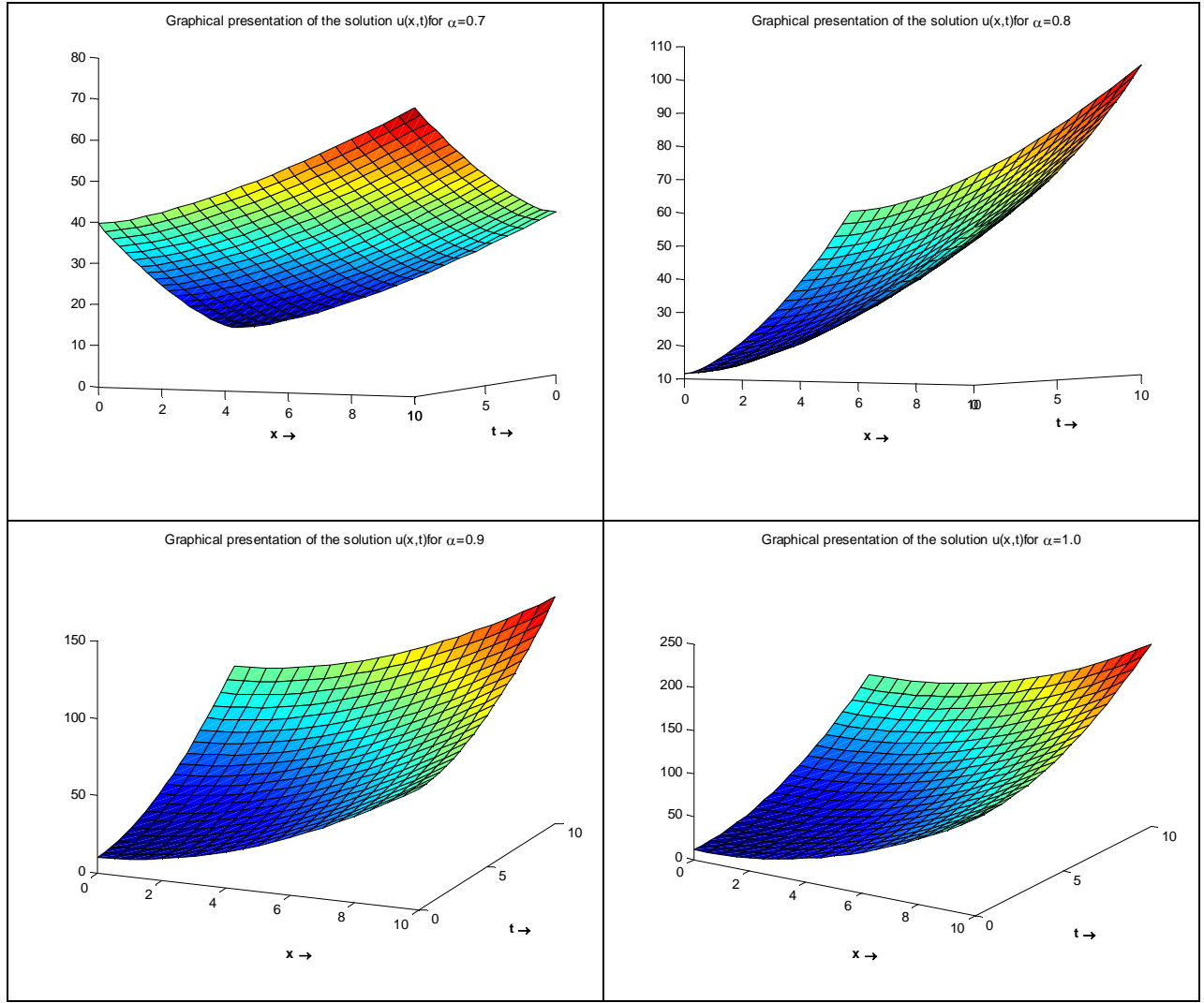

**Fig-1: Graphical representation of the solution (3.33) for** $\alpha = 0.7, 0.8, 0.9, 1.0$.

**Solution**: Considering $p = q = 1$ (as these constants are scale factors) in (3.31) we get the D'Alembert's solution of (3.32) in the following form

$$u(x,t) = \frac{1}{2}\left[ f\left(\left(\tfrac{1}{\Gamma(1+\alpha)} x^\alpha + \tfrac{1}{\Gamma(1+\alpha)} c^\alpha t^\alpha\right)\right) + f\left(\left(\tfrac{1}{\Gamma(1+\alpha)} x^\alpha - \tfrac{1}{\Gamma(1+\alpha)} c^\alpha t^\alpha\right)\right)\right] + \frac{1}{2c^\alpha}\int_{\left(\tfrac{1}{\Gamma(1+\alpha)} x^\alpha - \tfrac{1}{\Gamma(1+\alpha)} c^\alpha t^\alpha\right)}^{\left(\tfrac{1}{\Gamma(1+\alpha)} x^\alpha + \tfrac{1}{\Gamma(1+\alpha)} c^\alpha t^\alpha\right)} g(\xi)d\xi$$

Putting the initial conditions $f(x) = x^2$ and $g(x) = \sin x$ we get the following

$$\begin{aligned} u(x,t) &= \frac{1}{2}\left[\left(\tfrac{1}{\Gamma(1+\alpha)} x^\alpha + \tfrac{1}{\Gamma(1+\alpha)} c^\alpha t^\alpha\right)^2 + \left(\tfrac{1}{\Gamma(1+\alpha)} x^\alpha - \tfrac{1}{\Gamma(1+\alpha)} c^\alpha t^\alpha\right)^2\right] \\ &\quad + \frac{1}{2c^\alpha}\int_{\left(\tfrac{1}{\Gamma(1+\alpha)} x^\alpha - \tfrac{1}{\Gamma(1+\alpha)} c^\alpha t^\alpha\right)}^{\left(\tfrac{1}{\Gamma(1+\alpha)} x^\alpha + \tfrac{1}{\Gamma(1+\alpha)} c^\alpha t^\alpha\right)} \sin(\xi)d\xi \\ &= \left(\tfrac{1}{\Gamma(1+\alpha)} x^\alpha\right)^2 + \left(\tfrac{1}{\Gamma(1+\alpha)} c^\alpha t^\alpha\right)^2 + \frac{1}{c^\alpha}\cos\left(\tfrac{1}{\Gamma(1+\alpha)} x^\alpha\right)\cos\left(\tfrac{1}{\Gamma(1+\alpha)} c^\alpha t^\alpha\right) \end{aligned} \quad (3.33)$$

Thus we obtain solution of D'Alembert's of the fractional wave equation.

If we consider the initial condition as $f(x)=0$ and $g(x)=\sin x$ then we get the solution as

$$u(x,t)=\frac{1}{2c^{\alpha}}\int_{\left(\frac{1}{\Gamma(1+\alpha)}x^{\alpha}-\frac{1}{\Gamma(1+\alpha)}c^{\alpha}t^{\alpha}\right)}^{\left(\frac{1}{\Gamma(1+\alpha)}x^{\alpha}+\frac{1}{\Gamma(1+\alpha)}c^{\alpha}t^{\alpha}\right)}\sin(\xi)d\xi = \frac{1}{c^{\alpha}}\cos\left(\tfrac{1}{\Gamma(1+\alpha)}x^{\alpha}\right)\cos\left(\tfrac{1}{\Gamma(1+\alpha)}c^{\alpha}t^{\alpha}\right) \qquad (3.34)$$

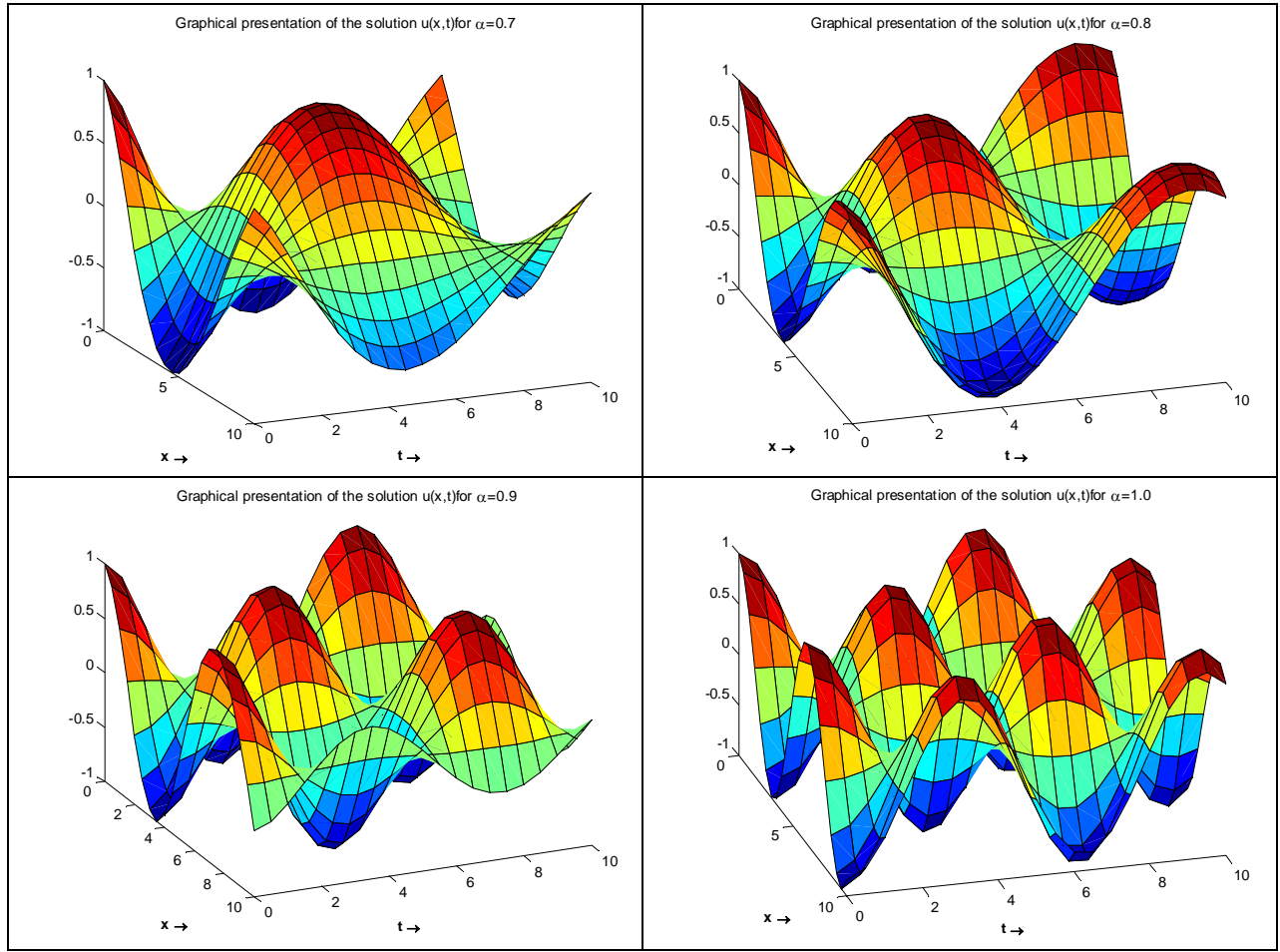

**Fig-2: Graphical representation of the solution (3.34) for $\alpha = 0.7, 0.8, 0.9, 1.0$**

From the graphical presentations in figure-1 and 2 it is clear that the solution depends on the order of fractional derivative, with the increase of order of fractional derivative the solution pattern changes.

## 5.0 Conclusions

In this paper we obtained the solution of $\alpha-$ and $2\alpha-$ order fractional wave equation; composed via Jumarie type fractional derivative using complex fractional transformation. We have demonstrated that using this fractional complex transformation method, we obtained the solutions easily as compared to fractional method of characteristics. The solution of $2\alpha-$ order fractional wave equation is like the D'Alembert's solution of classical second order wave equation. We show that it is having travelling wave components but are transformed into scaled coordinates $\sim x^\alpha$ and $\sim t^\alpha$. When the fractional order $\alpha$ is unity, we get the classical solution of the wave equation.